\documentclass[]{article}
\usepackage{amsfonts}
\usepackage{theorem}
\newtheorem{definition}{Definition}
\newtheorem{theorem}{Theorem}
\newtheorem{corollary}{Corollary}
\newtheorem{lemma}{Lemma}

\usepackage[unicode,colorlinks,plainpages=false]{hyperref}
\newcommand{\openbox}{$\begin{array}{c}
\hspace*{-0.55em}\sqcap \hspace*{-0.60em}\\[-0.4em] \hline
\multicolumn{1}{c}{\hspace*{-0.60em}}\\[-0.8em]
\end{array}
$}

\usepackage{enumerate}

\begin{document}

\centerline{{\bf On Commutative Monoid Congruences}
\footnote{Research supported by the Hungarian NFSR grant No T042481 and
No T043034}}
\centerline{\bf of Semigroups}

\bigskip

\bigskip

\centerline{\bf Attila Nagy}

\bigskip

\begin{abstract}
A subset $A$ of a semigroup $S$ is called a medial subset
of $S$ if $xaby\in A$ implies $xbay\in A$ for every $a,b,x,y\in S$.
By the separator of a subset $A$ of a semigroup $S$, we mean the set
of all elements $x$ of $S$ which satisfy the conditions
$xA\subseteq A$, $Ax\subseteq A$, $x\overline A\subseteq \overline A$,
$\overline Ax\subseteq \overline A$, where $\overline A$ denotes the
complement of $A$ in $S$.
In this paper we show that
if $\{A_i,\ i\in I\}$ is a family of medial subsets of a
semigroup $S$ such that $A=\cap _{i\in I}Sep(A_i)$ is not empty then
$P_{\{A_i,\ i\in I\} }$ defined by
$(a,b)\in P_{\{ A_i, i\in I\}}$ ($a,b\in S$)
if and only if, for every $i\in I$ and $x,y\in S$,
$xay\in A_i\ \Leftrightarrow \ xby\in A_i$ is a commutative monoid
congruence of $S$ such that $A$ is the identity element of
$S/P_{\{A_i,\ i\in I\} }$. Conversely, every commutative monoid congruence
of a semigroup can be so constructed.
We also show that if $S$ is a permutative semigroup then
the monoid congruences of $S$ are exactly the congruences
$P_{\{A_i,\ i\in I\} }$ defined
for arbitrary family $\{A_i,\ i\in I\}$ of arbitrary
subsets of $S$ satisfying $\cap _{i\in I}Sep(A_i)\neq \emptyset$.
\end{abstract}

\bigskip

By the idealizer $Id(A)$ of a subset $A$ of a semigroup $S$ we mean the
set of all elements $x$ of $S$ which satisfy the conditions $xA\subseteq A$,
$Ax\subseteq A$.
Denoting the complement of $A$ in $S$ by $\overline A$,
the subset $Sep(A)=Id(A)\cap Id(\overline A)$ of $S$ is called
the separator of $A$ (\cite{2}). In other words, the separator of $A$ is the set
of all elements $x$ of $S$ which satisfy the conditions $xA\subseteq A$,
$Ax\subseteq A$, $x\overline A\subseteq \overline A$, $\overline Ax\subseteq
\overline A$.

\begin{lemma}\label{lm1} (\cite{2}) For any subset $A$ of a semigroup $S$,
$Sep(A)$ is either empty or a subsemigroup of $S$.\hfill\openbox
\end{lemma}

\begin{lemma}\label{lm2} (\cite{2}) If $A$ is a subset of a semigroup such that
$Sep(A)\neq \emptyset$ then either $Sep(A)\subseteq A$ or
$Sep(A)\subseteq \overline A$.\hfill\openbox
\end{lemma}

A subset $U$ of a semigroup $S$ is said to be a left (right) unitary subset
of $S$ if $a,ab\in U$ ($a,ba\in U$) implies $b\in U$ for every $a,b\in S$.
The subset $U$ is called a unitary subset of $S$ if it is both left and right
unitary in $S$.

\begin{lemma}\label{lm3} (\cite{2})
A subsemigroup $A$ of a semigroup $S$ is unitary in $S$ if and only if
$A=Sep(A)$.\hfill\openbox
\end{lemma}

\newpage

Let $\{ A_i, i\in I\}$ be a family of non-empty subsets of a semigroup $S$.
It is easy to see that the relation $P_{\{ A_i, i\in I\}}$ on $S$ defined by
$(a,b)\in P_{\{ A_i, i\in I\}}$ ($a,b\in S$)
if and only if, for every $i\in I$ and $x,y\in S$,
$xay\in A_i\ \Leftrightarrow \ xby\in A_i$ is a congruence of $S$.

\begin{definition}\label{df1} A subset $A$ of a semigroup $S$ will be called a medial subset
of $S$ if $xaby\in A$ if and only if $xbay\in A$ for every $a,b,x,y\in S$.
\end{definition}

\begin{theorem}\label{th1} Let $\{A_i,\ i\in I\}$ be a family of medial subsets of a
semigroup $S$ such that $A=\cap _{i\in I}Sep(A_i)$ is not empty. Then
$P_{\{A_i,\ i\in I\} }$ is a commutative monoid
congruence of $S$ such that $A$ is the identity element of $S/P_{\{A_i,\ i\in I\} }$.
Conversely, every commutative monoid congruence
of a semigroup can be so constructed.
\end{theorem}

\noindent

{\bf Proof}. Let $\{A_i,\ i\in I\}$ be a family of medial subsets of a
semigroup $S$ such that $A=\cap _{i\in I}Sep(A_i)$ is not empty.
As $xaby\in A_i$ iff $xbay\in A_i$ for every $a,b,x,y\in S$ and $i\in I$,
we get that $P_{\{A_i,\ i\in I\} }$ is a commutative congruence on $S$. Let $a$ and $b$
be arbitrary elements of $S$ such that $a\in A_i$ and $b\notin A_i$ for some $i\in I$.
Then, for every $g,h\in A$, we have $gah\in A_i$ and $gbh\notin A_i$ and so
$(a,b)\notin P_{\{A_i,\ i\in I\} }$. Thus $A_i$ is
a union of $P_{\{A_i,\ i\in I\} }$-classes for every $i\in I$. Let $a,b\in A$
be arbitrary elements. Assume $xay\in A_i$ for some $i\in I$ and $x,y\in S$.
Since $b\in Sep(A_i)$, we get $xayb\in A_i$. As $A_i$ is a medial subset of $S$,
we get $xyab\in A_i$ and so $xy\in A_i$, because $ab\in Sep(A_i)$. Then
$xyba\in A_i$, $xbya\in A_i$ and $xby\in A_i$, because $ba\in Sep(A_i)$, $A_i$ is
medial and $a\in Sep(A_i)$. Thus
$(a,b)\in P_{\{A_i,\ i\in I\} }$. Let $a\in A$ and $b\notin A$ be arbitrary
elements. Then there is an index $j\in I$ such that $b\notin Sep(A_j)$.
We have four cases: $bA_j\not \subseteq A_j$, $A_jb\not \subseteq A_j$,
$b\overline A_j\not \subseteq \overline A_j$,
$\overline A_jb\not \subseteq \overline A_j$.
In case $bA_j\not \subseteq A_j$, there is an element $c\in A_j$ such that
$bc\notin A_j$ and so $abc\notin A_j$. As $aac\in A_j$, we get $(a,b)\notin
P_{\{A_i,\ i\in I\} }$. We get the same result in the other
three cases. Thus $A$ is a $P_{\{A_i,\ i\in I\} }$-class.
Let $a\in A$ and $s\in S$ be arbitrary elements. Then, for every $x,y\in S$,
$xsay\in A_i$ iff $xsaya\in A_i$ iff $xsyaa\in A_i$ iff
$xsy\in A_i$. Thus $(sa,s)\in P_{\{A_i,\ i\in I\} }$. We can prove, in a
similar way, that $(as,s)\in P_{\{A_i,\ i\in I\} }$. Hence $A$ is
the identity element in the factor semigroup $S/P_{\{A_i,\ i\in I\} }$.
Hence $P_{\{A_i,\ i\in I\} }$ is a commutative monoid congruence of $S$.

Conversely, let $\sigma$ be a commutative monoid congruence of a semigroup
$S$. Let $A$ denote the $\sigma$-class which is the identity element of
$S/\sigma$. Let $\{ A_i,\ i\in I\}$ denote the family of all $\sigma$-classes
of $S$. It is obvious that $A_i,\ i\in I$ are medial subsets of $S$. Let
$a\in A$ be an arbitrary element. As $aA_i\subseteq A_i$ and
$A_ia\subseteq A_i$ for every $i\in I$, we get $a\in \cap_{i\in I}Sep(A_i)$.
Hence $A\subseteq \cap_{i\in I}Sep(A_i)$. Assume that there is an element
$b$ of $S$ such that $b\in \cap_{i\in I}Sep(A_i)$ and $b\notin A$. Then there is
an index $j\in I$ such that $b\in A_j\neq A$ and so
$A_j\cap Sep(A_j)\neq \emptyset$. Then, by Lemma~\ref{lm2}, $Sep(A_j)\subseteq A_j$ which
implies $A\subseteq A_j$ which is impossible. Hence $A=\cap_{i\in I}Sep(A_i)$.
In the first part of the proof, it was proved that $A_i,\ i\in I$ are unions of
$P_{\{A_i,\ i\in I\} }$-classes. Hence $P_{\{A_i,\ i\in I\} }\subseteq \sigma$.
As every $A_i, i\in I$ is a $\sigma$-class, it is obvious that
$\sigma \subseteq P_{\{A_i,\ i\in I\} }$. Consequently
$\sigma =P_{\{A_i,\ i\in I\} }$.
\hfill{\openbox}

\bigskip
A subset $A$ of a semigroup $S$ is said to be a reflexive subset of $S$
if $ab\in A$ implies $ba\in A$ for every $a,b\in S$.

\medskip

\begin{corollary}\label{cr1}
For any medial subset $A$ of a semigroup $S$, $Sep(A)$ is either empty
or a reflexive unitary subsemigroup of $S$.
\end{corollary}

\noindent

{\bf Proof}. Let $A$ be a medial subset of a semigroup $S$ such that
$Sep(A)\neq \emptyset$. By Theorem~\ref{th1}, $P_A$ is a monoid congruence of $S$
such that $Sep(A)$ is the identity element of the factor semigroup $S/P_A$.
Then $Sep(Sep(A))=Sep(A)$ and, by Lemma~\ref{lm3}, $Sep(A)$ is a unitary subsemigroup
of $S$. As $P_A$ is a commutative congruence by Theorem~\ref{th1}, $Sep(A)$ is reflexive.
\hfill{\openbox}

\begin{definition} A semigroup $S$ is called a permutative semigroup (\cite{3}) if
it satisfies a non-trivial permutation identity, that is, there is a
positive integer $n\geq 2$ and a non-identity permutation $\sigma$ of
$\{ 1, 2,\dots ,n\}$ such that $S$ satisfies the identity $x_1x_2\dots x_n=
x_{\sigma (1)}x_{\sigma (2)}\dots x_{\sigma (n)}$.
\end{definition}

\medskip

It is obvious that every permutative monoid is commutative. Next, we
construct the monoid congruences of permutative semigroups.

\begin{lemma}\label{lm4} {\cite{4}} Let $S$ be a permutative semigroup. Then there exists a
positive integer $k$ such that, for every $u,v\in S^k$ and $x,y\in S$, we
have $uxyv=uyxv$.
\end{lemma}

\begin{theorem}\label{th2}
Let $\{A_i,\ i\in I\}$ be a family of subsets of a permutative
semigroup $S$ such that $A=\cap _{i\in I}Sep(A_i)$ is not empty. Then
$P_{\{A_i,\ i\in I\} }$ is a monoid congruence
of $S$ such that $A$ is the identity element of $S/P_{\{A_i,\ i\in I\} }$.
Conversely, every monoid congruence of a permutative semigroup can be
so constructed.
\end{theorem}

\noindent

{\bf Proof}. Let $S$ be a permutative semigroup. Then, by Lemma~\ref{lm4},
there is a positive integer $k$ such that, for every $u,v\in S^k$ and $x,y\in S$, we
have $uxyv=uyxv$. Let $X$ be a non-empty subset
of $S$ such that $Sep(X)\neq \emptyset$. Assume $uxyv\in X$ for some $u,v,x,y\in S$.
Then, for some $t\in Sep(X)$, we have $(t^{k-1}u)yx(vt^{k-1})=(t^{k-1}u)xy(vt^{k-1})\in X$
which implies $uyxv\in X$. Hence $X$ is a medial subset of $S$.
Assume that $\{ A_i,i\in I\}$
is a family of subsets of $S$ such that  $A=\cap _{i\in I}Sep(A_i)$ is not
empty. Then, by the above, every $A_i$ is a medial subset of $S$ and so, by
Theorem~\ref{th1}, $P_{\{A_i,\ i\in I\} }$ is a (commutative) monoid
congruence of $S$ such that $A$ is the identity element of
$S/P_{\{A_i,\ i\in I\} }$. The converse follows from Theorem~\ref{th1}.
\hfill{\openbox}

\begin{corollary}\label{cr2}
For any subset $A$ of a permutative semigroup $S$, $Sep(A)$ is either empty
or a reflexive unitary subsemigroup of $S$.
\end{corollary}

\noindent

{\bf Proof}. As a subset $A$ of a permutative semigroup $S$ with
$Sep(A)\neq \emptyset$ is medial, the assertion follows from Corollary~\ref{cr1}.
\hfill{\openbox}

\medskip

\noindent
Attila Nagy

\noindent
Department of Algebra

\noindent
Mathematical Institute

\noindent
Budapest University of Technology and Economics

\noindent
e-mail: nagyat@math.bme.hu

\bigskip


\bigskip

\begin{thebibliography}{4}

\bibitem{1} Clifford, A.H. and G.B. Preston, {\cal The Algebraic Theory
of Semigroups}, Amer. Math. Soc., Providence, R.I., I(1961), II(1967)

\bibitem{2} Nagy, A., {\cal The separator of a subset of a semigroup},
Publicationes Mathematicae, Tom. 27., Fasc. 1-2.(1980), 25-30

\bibitem{3} Nordahl, T.E., {\cal On permutative semigroup algebras},
Algebra Universalis, 25(1988), 322-333

\bibitem{4} Putcha, M.S. and A. Yaqub, {\cal Semigroups satisfying permutation
identities}, Semigroup Forum, 3(1971), 68-73


\end{thebibliography}
\end{document}